\documentclass{amsart}%
\usepackage{amsfonts}
\usepackage{amsmath}
\usepackage{amssymb}
\usepackage[dvips]{graphicx}%
\newtheorem{theorem}{Theorem}[section]
\theoremstyle{plain}
\newtheorem{corollary}[theorem]{Corollary}
\newtheorem{lemma}[theorem]{Lemma}
\newtheorem{proposition}[theorem]{Proposition}
\theoremstyle{definition}

\newtheorem{definition}[theorem]{Definition}
\newtheorem{definitions}[theorem]{Definitions}
\newtheorem{example}[theorem]{Example}

\newtheorem{remark}[theorem]{Remark}
\newtheorem{remarks}[theorem]{Remarks}
\numberwithin{equation}{section}
\renewcommand{\theenumi}{\alph{enumi}}
\newcounter{enump}
\renewcommand{\theenump}{\roman{enump}}
\newlength{\qedskip}
\newlength{\qedadjust}
\def\hooklongrightarrow{\mathrel{\lhook\joinrel\longrightarrow}} 
\begin{document}
\title[Frames in Hilbert space and reproducing kernels]{Frame analysis and approximation in\\reproducing kernel Hilbert spaces}
\author{Palle E. T. Jorgensen}
\address{Department of Mathematics, The University of Iowa, 14 MacLean Hall, Iowa City,
IA 52242-1419, U.S.A.}
\email{jorgen@math.uiowa.edu}
\urladdr{http://www.math.uiowa.edu/\symbol{126}jorgen}
\thanks{This material is based upon work supported by the U.S. National Science
Foundation under grants DMS-0139473 (FRG) and DMS-0457581.}
\subjclass[2000]{33C50, 42C15, 46E22, 47B32}
\keywords{tight frames, Parseval frames, Bessel sequences, Hilbert space, reproducing
kernel, Shannon, polar decomposition}

\begin{abstract}
We consider frames $F$ in a given Hilbert space, and we show that every $F$
may be obtained in a constructive way from a reproducing kernel and an
orthonormal basis in an ambient Hilbert space. The construction is
operator-theoretic, building on a geometric formula for the analysis operator
defined from $F$. Our focus is on the infinite-dimensional case where \emph{a
priori} estimates play a central role, and we extend a number of results which
are known so far only in finite dimensions. We further show how this approach
may be used both in constructing useful frames in analysis and in
applications, and in understanding their geometry and their symmetries.

\end{abstract}
\maketitle

\section{\label{Int}Introduction}

Frames are redundant bases which turn out in certain applications to be more
flexible than the better known orthonormal bases (ONBs) in Hilbert space. The
frames allow for more symmetries than ONBs do, especially in the context of
signal analysis, and of wavelet constructions; see, e.g.,
\cite{CoDa93,BDP05,Dut06a}. Since frame bases (although containing
redundancies) still allow for efficient algorithms, they have found many
applications, even in finite dimensions; see, for example,
\cite{BeFi03,CaCh03,Chr03,Eld02,FJKO05,VaWa05}.

As is well known, when a vector $f$ in a Hilbert space $\mathcal{H}$ is
expanded in an orthonormal basis $B$, there is then automatically an
associated Parseval identity. In physical terms, this identity typically
reflects a \emph{stability} feature of a decomposition based on the chosen ONB
$B$. Specifically, Parseval's identity reflects a conserved quantity for a
problem at hand, for example, energy conservation in quantum mechanics.

The theory of frames (see Definitions \ref{DefMor.1}) begins with the
observation that there are useful vector systems which are in fact not ONBs
but for which a Parseval formula still holds. In fact, in applications it is
important to go beyond ONBs. While this viewpoint originated in signal
processing (in connection with frequency bands, aliasing, and filters), the
subject of frames appears now to be of independent interest in mathematics.

On occasion, we may have a system of vectors $S$ in $\mathcal{H}$ for which
Parseval's identity is still satisfied, but such that a generalized Parseval's
identity might only hold up to a fixed constant $c$ of scale. (For example, in
sampling theory, a scale might be introduced as a result of
``oversampling''.) In this case, we say that the constant $c$ scales the
expansion. Suppose a system of vectors $S$ in a given Hilbert space
$\mathcal{H}$ allows for an expansion, or decomposition of every $f$ in
$\mathcal{H}$, but the analogue of Parseval's identity holds only up to a
fixed constant $c$ of scale. In that case, we say that $S$ is a \emph{tight
frame} with frame constant $c$. So the special case $c = 1$ is the case of a
Parseval frame. For precise definitions of these terms, we refer to Section
\ref{Pre} below, or to the book literature, e.g., \cite{Chr03}.

Aside from applications, at least three of the other motivations for frame
theory come from: (1)~wavelets, e.g., \cite{CoDa93}, \cite{BJMP05}, and
\cite{BJMP06}; (2)~from non-harmonic Fourier expansions \cite{DuSc52}; and
(3)~from computations with polynomials in several variables, and their
generalized orthogonality relations \cite{DuXu01}.

While frames already have impressive uses in signal processing (see, e.g.,
\linebreak
\cite{ALTW04,Chr99}), they have recently \cite{CCLV05,CKL04} been shown to be
central in our understanding of a fundamental question in operator algebras,
the Kadison--Singer conjecture. We refer the reader to \cite{CFTW06} for
up-to-date research, and to \cite{Chr99,KaRi97,Nel57,Nel59} for background.

In all these cases, the authors work with recursive algorithms, and the issue
of \emph{stability} plays a crucial role. Stability, however, may obtain in
situations that are much more general than the context of traditional ONBs, or
even tight frames. In fact, stability may apply even when we have only \emph{a
priori} estimates, as opposed to identities: for example, when the scaled
version of Parseval's identity is replaced with a pair of estimates, a fixed
lower bound and an upper bound; see (\ref{eqMor.1}) below. If such bounds
exist, they are called lower and upper frame bounds.

If a system $S$ of vectors in a Hilbert space $\mathcal{H}$ satisfies such a
pair of \emph{a priori} estimates, we say that $S$ is simply a \emph{frame}.
And if such an estimate holds only with an \emph{a priori} upper bound, we say
that $S$ is a \emph{Bessel sequence}. It is known (see, e.g., \cite{AkGl93})
that for a fixed Hilbert space $\mathcal{H}$, the various classes of frames
$S$ in $\mathcal{H}$ may be obtained from some ambient Hilbert space
$\mathcal{K}$ and an orthonormal basis $B$ in $\mathcal{K}$, i.e., when the
pair $(S, \mathcal{H})$ is given, there are choices of $\mathcal{K}$ such that
the frame $S$ may be obtained from applying a certain bounded operator $T$ to
a suitable ONB $B$ in $\mathcal{K}$. Passing from the given structure in
$\mathcal{H}$ to the ambient Hilbert space is called \emph{dilation} in
operator theory. The properties of the operator $T$ which does the job depend
on the particular frame in question. For example, if $S$ is a Parseval frame,
then $T$ will be a projection of the ambient Hilbert space $\mathcal{K}$ onto
$\mathcal{H}$. But this operator-theoretic approach to frame theory has been
hampered by the fact that the ambient Hilbert space is often an elusive
abstraction. Starting with a frame $S$ in a fixed Hilbert space $\mathcal{H}$,
then by dilation, or extension, we pass to an ambient Hilbert space
$\mathcal{K}$. In this paper we make concrete the selection of the ``magic''
operator $T \colon\mathcal{K} \to\mathcal{H}$ which maps an ONB in
$\mathcal{K}$ onto $S$. While existence is already known, the building of a
dilation system $(\mathcal{K}, T, \text{ONB})$ is often rather
non-constructive, and the various methods for getting $\mathcal{K}$ are
fraught with choices that are not unique.

Nonetheless, it was shown recently \cite{Dut04c,Dut06a} that when the dilation
approach is applied to Parseval frames of wavelets in $\mathcal{H}%
=L^{2}(\mathbb{R})$, i.e., to wavelet bases which are not ONBs, then the
ambient Hilbert space $\mathcal{K}$ can be made completely explicit, and the
constructions are algorithmic. Moreover, the \textquotedblleft
inflated\textquotedblright\ ONB in $\mathcal{K}$ then takes the form of a
traditional ONB-wavelet basis, a so-called \textquotedblleft super-%
wavelet\textquotedblright. For details, see \cite{BDP05,Dut04c}, and also the
papers \cite{KoLa04,BJMP05,BJMP06}.

It is the purpose of the present paper to show that the techniques which work
well in this restricted context, \textquotedblleft super-%
wavelets\textquotedblright\ and redundant wavelet frames, apply to a more
general and geometric context, one which is motivated in turn by extension
principles in probability theory; see, e.g., \cite{PaSc72}, \cite{Dut04b}, and
\cite{Jor06a}.

A key idea in our present approach is the use of reproducing Hilbert spaces,
and their reproducing kernels in the sense of \cite{Aro50}. See also
\cite{Nel57,Nel59} for an attractive formulation. We show that for every
Hilbert space $\mathcal{H}$, and every frame $S$ in $\mathcal{H}$ (even if $S$
is merely a Bessel sequence), there is a way of constructing the ambient
Hilbert space $\mathcal{K}$ in such a way that the operator $T$ has a concrete
reproducing kernel.

Finally, we mention that a recent paper \cite{VaWa05} serves as a second
motivation for our work; in fact \cite{VaWa05} contains finite-dimensional
cases of two of our present theorems. These results in \cite{VaWa05} are
Theorems 2.9 and 6.4 in that paper: The results in \cite{VaWa05} are concerned
with symmetries of tight frames, and with associated families of unitary
representations of the symmetry groups. It turns out that this approach to
symmetry is natural in the context of operator theory; see, e.g.,
\cite{PaSc72}, \cite{Dut04b}, \cite{Dut06a}.

\section{\label{Pre}Preliminary notions and definitions}

Let $S$ be a countable set, finite or infinite, and let $\mathcal{H}$ be a
complex or real Hilbert space. We shall be interested in a class of spanning
families of vectors $\left(  \mathbf{v}\left(  s\right)  \right)  $ in
$\mathcal{H}$ indexed by points $s\in S$. Their properties will be defined
precisely below, and the families are termed \emph{frames}. The simplest
instance of this is when $\mathcal{H}=\ell^{2}\left(  S\right)  ={}$the
Hilbert space of all square-summable sequences, i.e., all $f\colon
S\rightarrow\mathbb{C}$ such that $\sum_{s\in S}\left\vert f\left(  s\right)
\right\vert ^{2}<\infty$. In that case, set%
\begin{equation}
\left\langle \,f_{1}\mid f_{2}\,\right\rangle :=\sum_{s\in S}\,\overline
{f_{1}\left(  s\right)  }\,f_{2}\left(  s\right)  \label{eqPre.1}%
\end{equation}
for all $f_{1},f_{2}\in\ell^{2}\left(  S\right)  $.

It is then immediate that the delta functions $\left\{  \,\delta_{s}\mid s\in
S\,\right\}  $ given by%
\begin{equation}
\delta_{s}\left(  t\right)  =%
\begin{cases}
1, & t=s,\\
0, & t\in S\setminus\left\{  s\right\}  ,
\end{cases}
\label{eqPre.2}%
\end{equation}
form an \emph{orthonormal basis} (ONB) for $\mathcal{H}$, i.e., that%
\begin{equation}
\left\langle \,\delta_{s_{1}}\mid\delta_{s_{2}}\,\right\rangle =%
\begin{cases}
1 & \text{if }s_{1}=s_{2}\text{ in }S,\\
0 & \text{if }s_{1}\neq s_{2},
\end{cases}
\label{eqPre.3}%
\end{equation}
and that this is a maximal orthonormal family in $\mathcal{H}$. Moreover,%
\begin{equation}
f=\sum_{s\in S}f\left(  s\right)  \delta_{s}\text{\qquad for all }f\in\ell
^{2}\left(  S\right)  . \label{eqPre.decomp}%
\end{equation}

It also is immediate from (\ref{eqPre.1}) that Parseval's formula%
\begin{equation}
\left\Vert f\right\Vert ^{2}=\sum_{s\in S}\left\vert \left\langle \,\delta
_{s}\mid f\,\right\rangle \right\vert ^{2} \label{eqPre.4}%
\end{equation}
holds for all $f\in\ell^{2}\left(  S\right)  $.

We shall consider pairs $\left(  S,\mathcal{H}\right)  $ and indexed families
\begin{equation}
\left\{  \,\mathbf{v}\left(  s\right)  \mid s\in S\,\right\}  \subset
\mathcal{H} \label{eqPre.5}%
\end{equation}
such that for some $c\in\mathbb{R}_{+}$, the identity%
\begin{equation}
\left\Vert f\right\Vert ^{2}=c\sum_{s\in S}\left\vert \left\langle
\,\mathbf{v}\left(  s\right)  \mid f\,\right\rangle \right\vert ^{2}
\label{eqPre.6}%
\end{equation}
holds for all $f\in\mathcal{H}$.

When a Hilbert space $\mathcal{H}$ is given, our main result states that
solutions to (\ref{eqPre.6}) exist if and only if $\mathcal{H}$ is
isometrically embedded in $\ell^{2}\left(  S\right)  $. But we further
characterize these embeddings, and we use this in understanding the geometry
of tight frames (details below).

\begin{definition}
\label{DefPre.1}Let $\left(  S,\mathcal{H},c,\left(  \mathbf{v}\left(
s\right)  \right)  _{s\in S}\right)  $ be as above.
We shall say that this system constitutes a \emph{tight frame} with frame
constant $c$ if \textup{(\ref{eqPre.6})} holds.
\end{definition}

(Note that if $\left(  \mathbf{v}\left(  s\right)  \right)  _{s\in S}$
satisfies (\ref{eqPre.6}), then the scaled system $\left(  \sqrt
{c}\,\mathbf{v}\left(  s\right)  \right)  _{s\in S}$ has the property with
frame constant one.)

\begin{example}
\label{ExaPre.2}\textup{($S$ finite.)} Let $\mathcal{H}$ be the
two-dimensional real Hilbert space, and let $n\geq3$. Set $S:=\left\{
1,2,\dots,n\right\}  =:S_{n}$, and%
\begin{equation}
\mathbf{v}\left(  s\right)  :=%
\begin{pmatrix}
\displaystyle\cos\left(  \frac{2\pi s}{n}\right)  _{\mathstrut}\\
\displaystyle\sin\left(  \frac{2\pi s}{n}\right)  ^{\mathstrut}%
\end{pmatrix}
,\qquad s\in S. \label{eqPre.7}%
\end{equation}
Then it is easy to see that this constitutes a tight frame with frame constant
$c=\frac{2}{n}$. Examples are presented in Figures \textup{\ref{FigExaPre.2.1}%
} and \textup{\ref{FigExaPre.2.2}}.\begin{figure}[ptb]
\setlength{\unitlength}{108bp} \begin{picture}(3.1,1)
\put(0,0){\includegraphics[bb=88 4 196 112,width=\unitlength]{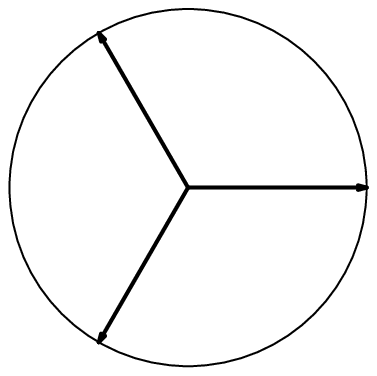}}
\put(2,0){\includegraphics[bb=88 4 196 112,width=\unitlength]{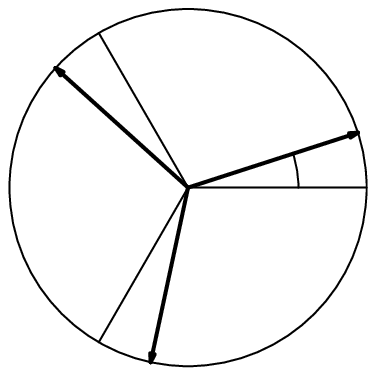}}
\put(2.80525,0.54809){\makebox(0,0)[l]{$\,\theta$}}
\end{picture}
\caption{Two illustrations for $n=3$}%
\label{FigExaPre.2.1}%
\end{figure}\begin{figure}[ptbptb]
\setlength{\unitlength}{108bp} \begin{picture}(3.1,1)
\put(0,0){\includegraphics[bb=88 4 196 112,width=\unitlength]{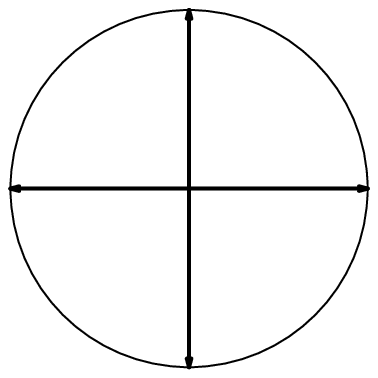}}
\put(2,0){\includegraphics[bb=88 4 196 112,width=\unitlength]{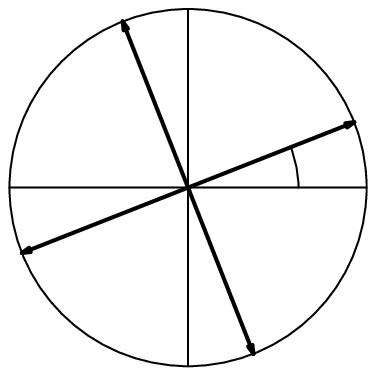}}
\put(2.80360,0.55760){\makebox(0,0)[l]{$\,\theta$}}
\end{picture}
\caption{Two illustrations for $n=4$}%
\label{FigExaPre.2.2}%
\end{figure}
\end{example}

\begin{definition}
\label{DefPre.3}Let $\mathcal{H}$ and $\mathcal{K}$ be Hilbert spaces over
$\mathbb{C}$ or $\mathbb{R}$, and let $V\colon\mathcal{H}\rightarrow
\mathcal{K}$ be a linear mapping. We say that $V$ is an \emph{isometry}, and
that $\mathcal{H}$ is \emph{isometrically embedded} in $\mathcal{K}$
\textup{(}via $V$\/\textup{)} if%
\begin{equation}
\left\Vert Vf\right\Vert _{\mathcal{K}}=\left\Vert f\right\Vert _{\mathcal{H}%
},\qquad f\in\mathcal{H}. \label{eqPre.8}%
\end{equation}

\end{definition}

Given a linear operator $V\colon\mathcal{H}\rightarrow\mathcal{K}$, we then
denote the adjoint operator $V^{\ast}\colon\mathcal{K}\rightarrow\mathcal{H}$.
It is easy to see that $V$ is \emph{isometric} iff $V^{\ast}V=I_{\mathcal{H}}$,
where $I_{\mathcal{H}}$ denotes the identity operator in $\mathcal{H}$.
Moreover, if $V$ is isometric then $P=P_{V}=VV^{\ast}\colon\mathcal{K}%
\rightarrow\mathcal{K}$ is a projection, i.e.,
\begin{equation}
P=P^{\ast}=P^{2} \label{eqPre.9}%
\end{equation}
holds, and the subspace%
\begin{equation}
P\mathcal{K}\subset\mathcal{K} \label{eqPre.10}%
\end{equation}
may be identified with $\mathcal{H}$ via the isometric embedding.

We state our next result only in the case of frame constant $c=1$, but as
noted it easily generalizes.

\begin{theorem}
\label{ThmPre.4}Let $S$ be a countable set, and let $\mathcal{H}$ be a Hilbert
space over $\mathbb{C}$ \textup{(}or $\mathbb{R}$\textup{)}.

Then the following two conditions are equivalent:\renewcommand{\theenumi}{\roman{enumi}}

\begin{enumerate}
\item \label{ThmPre.4(1)}There is a tight frame $\left\{  \,\mathbf{v}\left(
s\right)  \mid s\in S\,\right\}  \subset\mathcal{H}$ with frame constant $c=1$;

\item \label{ThmPre.4(2)}$\mathcal{H}$ is isometrically embedded \textup{(}as
a closed subspace\/\textup{)} in $\ell^{2}\left(  S\right)  $.
\end{enumerate}
\end{theorem}

\begin{proof}
The details of proof will follow in the next section.
\end{proof}

\begin{remark}
\label{RemParseval}Tight frames with frame constant equal to one are called
\emph{Parseval frames}.
\end{remark}

\refstepcounter{theorem} \refstepcounter{theorem} \refstepcounter{theorem}

\section{\label{Pro}Proof of Theorem \ref{ThmPre.4}}

Let the data be as specified in the theorem: A pair $\left(  S,\mathcal{H}%
\right)  $ is given where $S$ is a set and $\mathcal{H}$ is a Hilbert space.
The conclusion is that (\refstepcounter{enump}\label{ThmPre.4proof1}%
\ref{ThmPre.4proof1})~the $S$-tight-frame property for $\mathcal{H}$ is
equivalent to (\refstepcounter{enump}\label{ThmPre.4proof2}%
\ref{ThmPre.4proof2})~the (isometric) embedding of $\mathcal{H}$ into
$\ell^{2}\left(  S\right)  $.

Assume (\ref{ThmPre.4proof1}), and let $\left\{  \,\mathbf{v}\left(  s\right)
\mid s\in S\,\right\}  \subset\mathcal{H}$ be the frame system which is
asserted. Now define%
\begin{equation}
V\colon\mathcal{H}\ni f\longmapsto\left(  \left\langle \,\mathbf{v}\left(
s\right)  \mid f\,\right\rangle _{\mathcal{H}}\right)  _{s\in S}.
\label{eqPro.1}%
\end{equation}
This operator is called the \emph{analysis operator} for the frame.

Equivalently, the analysis operator $V$ may be written in terms of the
canonical ONB $\left(  \delta_{s}\right)  _{s\in S}$ for $\ell^{2}\left(
S\right)  $ as follows:%
\begin{equation}
Vf=\sum_{s\in S}\left\langle \,\mathbf{v}\left(  s\right)  \mid
f\,\right\rangle \delta_{s}\qquad(\in\ell^{2}\left(  S\right)  ).
\label{eqProAnalysis}%
\end{equation}

Our first two observations are that the range of $V$ is contained in $\ell
^{2}\left(  S\right)  $, and that $V$ is isometric. But note that both these
conclusions are immediate consequences of identity (\ref{eqPre.6}) for the
special case $c=1$.

We formalize this in the following lemma. (The proof of converse implication
in Theorem \ref{ThmPre.4} will be resumed after the lemma.)

\begin{lemma}
\label{LemPro.half}Let $\left(  \mathbf{v}\left(  s\right)  \right)  _{s\in
S}$ be a system in the Hilbert space $\mathcal{H}$, and let $V$ be defined by
\textup{(\ref{eqPro.1})}. Then \textup{(\ref{eqPre.6})} holds if and only if
$\sqrt{c}\,V$ is isometric.
\end{lemma}

An easy calculation yields the formula for the adjoint operator:
\begin{equation}
V^{\ast}\left(  \left(  c_{s}\right)  _{s\in S}\right)  =\sum_{s\in S}%
c_{s}\mathbf{v}\left(  s\right)  ,\qquad\left(  c_{s}\right)  \in\ell
^{2}\left(  S\right)  . \label{eqPro.2}%
\end{equation}
Hence%
\begin{equation}
f=V^{\ast}Vf=\sum_{s\in S}\left\langle \,\mathbf{v}\left(  s\right)  \mid
f\,\right\rangle \mathbf{v}\left(  s\right)  \label{eqPro.3}%
\end{equation}
holds for all $f\in\mathcal{H}$.

Moreover, the projection operator $P:=VV^{\ast}$ is given by the formula%
\begin{equation}
\left(  P\left(  c_{s}\right)  \right)  _{t}=\sum_{s\in S}\left\langle
\,\mathbf{v}\left(  t\right)  \mid\mathbf{v}\left(  s\right)  \,\right\rangle
_{\mathcal{H}}c_{s}; \label{eqPro.4}%
\end{equation}
in other words, $P$ has a concrete matrix representation in the Hilbert space
$\ell^{2}\left(  S\right)  $. Specifically, $P$ is represented as
multiplication on column vectors $\left(  c_{s}\right)  _{s\in S}$, and the
matrix of $P$ is
\begin{equation}
P\left(  t,s\right)  =\left\langle \,\mathbf{v}\left(  t\right)
\mid\mathbf{v}\left(  s\right)  \,\right\rangle _{\mathcal{H}}.
\label{eqPro.5}%
\end{equation}

To prove the converse, assume (\ref{ThmPre.4proof2}).
Hence, there is an isometry $V\colon\mathcal{H}\rightarrow\ell^{2}\left(
S\right)  $. As we noted, this means that $P:=VV^{\ast}$ is a projection in
$\ell^{2}\left(  S\right)  $, and that%
\begin{equation}
P\ell^{2}\left(  S\right)  =\left\{  \,Vf\mid f\in\mathcal{H}\,\right\}  .
\label{eqPro.6}%
\end{equation}
Now observe that $\mathbf{v}\left(  s\right)  :=P\left(  \delta_{s}\right)  $,
where $\delta_{s}$ for each $s\in S$ is the delta function of (\ref{eqPre.2}).
We claim that (\ref{eqPre.6}) holds for $c=1$. To see this, let $f\in
\mathcal{H}$ be given. Then%
\begin{align*}
\left\Vert f\right\Vert ^{2}  &  =\left\Vert Vf\right\Vert ^{2}\underset
{\text{by (\ref{eqPre.4})}}{=}\sum_{s\in S}\left\vert \left\langle
\,\delta_{s}\mid Vf\,\right\rangle \right\vert ^{2}\\
&  =\sum_{s\in S}\left\vert \left\langle \,\delta_{s}\mid PVf\,\right\rangle
\right\vert ^{2}\\
&  =\sum_{s\in S}\left\vert \left\langle \,P\delta_{s}\mid Vf\,\right\rangle
\right\vert ^{2}\\
&  =\sum_{s\in S}\left\vert \left\langle \,\mathbf{v}\left(  s\right)  \mid
f\,\right\rangle \right\vert ^{2},
\end{align*}
which is the desired conclusion.\qed

\begin{remark}
\label{RemPro.Newpound}Even when no restricting assumptions are placed on a
given system of vectors $\left(  \mathbf{v}\left(  s\right)  \right)  _{s\in
S}$, the function $P$ from \textup{(\ref{eqPro.5})} is always positive
semidefinite; see Definition \textup{\ref{DefSym.1}}. Moreover, the following
converse implication holds \textup{(}Theorem \textup{\ref{ThmSym.2}}
below\/\textup{)}, known as the \emph{reconstruction principle\/}: Every positive semidefinite
function has the form \textup{(\ref{eqPro.5})} for some Hilbert space and an
associated system of vectors. As we show in Section \textup{\ref{Mor}}, when
this is specialized, we find that there is a graduated system of frame
properties which may or may not hold for a given system of vectors $\left(
\mathbf{v}\left(  s\right)  \right)  _{s\in S}$. Moreover we show that each of
these properties reflects a corresponding axiom for the associated function
\textup{(\ref{eqPro.5})}, Gramian, or correlation matrix. Our results in
Sections \textup{\ref{Sym}} and \textup{\ref{Mor}} below serve three purposes:
(1)~one is to use reproducing kernels \textup{\cite{Aro50,Nel57,PaSc72}} and
their spectral theory to study classes of frames and their symmetries;
(2)~another \textup{(}e.g., Theorems \textup{\ref{ThmSym.4}} and
\textup{\ref{ThmMor.2})} is to use an idea from operator theory to establish
results about frame deformations and their stability; and finally (3)~these results serve to tie
the operator theory to more current applications.
\end{remark}

\begin{corollary}
\label{CorPro.ter}For $i=1,2$, consider two systems $\left(  S_{i}%
,\mathcal{H}_{i}\right)  $ of sets and Hilbert spaces, and assume that
$\left(  \mathbf{v}_{i}\left(  s_{i}\right)  \right)  _{s_{i}\in S_{i}}$ is a
pair of tight frames with respective frame constants $c_{i}$. Then%
\begin{equation}
\mathbf{w}\left(  s_{1},s_{2}\right)  :=\mathbf{v}_{1}\left(  s_{1}\right)
\otimes\mathbf{v}_{2}\left(  s_{2}\right)  \label{eqPro.newpound}%
\end{equation}
defines a tight frame for the tensor-product Hilbert space $\mathcal{H}%
_{1}\otimes\mathcal{H}_{2}$ with frame constant $c=c_{1}c_{2}$.
\end{corollary}

\begin{proof}
The result follows immediately from the theorem and Lemma \ref{LemPro.half}.
To see this, note that if $V_{i}$, $i=1,2$, are the operators defined in
(\ref{eqPro.1}) above, then both%
\begin{equation}
\sqrt{c_{i}}\,V_{i}\colon\mathcal{H}_{i}\longrightarrow\ell^{2}\left(
S_{i}\right)  \label{eqCorPro.terproof1}%
\end{equation}
are isometric embeddings. It then follows that%
\begin{equation}
\sqrt{c_{1}}\,V_{1}\otimes\sqrt{c_{2}}\,V_{2}=\sqrt{c_{1}c_{2}}\,V_{1}\otimes
V_{2} \label{eqCorPro.terproof2}%
\end{equation}
is an isometric embedding of $\mathcal{H}_{1}\otimes\mathcal{H}_{2}$ into $\ell
^{2}\left(  S_{1}\times S_{2}\right)  $. From the fact that each $V_{i}$ is
defined from the system $\left(  \mathbf{v}_{i}\left(  s_{i}\right)  \right)
_{s_{i}\in S_{i}}$, it follows readily that the tensor operator $V_{1}\otimes
V_{2}$ is defined from the system in (\ref{eqPro.newpound}).

Specifically,
\begin{equation}
\mathbf{w}\left(  s_{1},s_{2}\right)  =\left(  V_{1}\otimes V_{2}\right)
^{\ast}\left(  \delta_{s_{1}}\otimes\delta_{s_{2}}\right)  =V_{1}^{\ast}%
\delta_{s_{1}}\otimes V_{2}^{\ast}\delta_{s_{2}}=\mathbf{v}_{1}\left(
s_{1}\right)  \otimes\mathbf{v}_{2}\left(  s_{2}\right)  ,
\label{eqCorPro.terproof3}%
\end{equation}
which is the desired result.
\end{proof}

\begin{definition}
\label{DefPro.1}Let $\left(  S,\mathcal{H}\right)  $ be a pair as above: $S$
is a set and $\mathcal{H}$ is a Hilbert space. Then the matrix%
\begin{equation}
k\left(  t,s\right)  :=\left\langle \,\mathbf{v}\left(  t\right)
\mid\mathbf{v}\left(  s\right)  \,\right\rangle _{\mathcal{H}},\qquad s,t\in
S, \label{eqPro.7}%
\end{equation}
is called the \emph{Gramian}, or the \emph{Gram matrix}.
\end{definition}

\begin{corollary}
\label{CorPro.2}Let the pair $\left(  S,\mathcal{H}\right)  $ be as in the
statement of Theorem \textup{\ref{ThmPre.4}}. Let $\left\{  \,\mathbf{v}%
\left(  s\right)  \mid s\in S\,\right\}  $ be a system of vectors in
$\mathcal{H}$ and let $k\left(  t,s\right)  =\left\langle \,\mathbf{v}\left(
t\right)  \mid\mathbf{v}\left(  s\right)  \,\right\rangle $ be the
corresponding Gram matrix. Then $\left(  \mathbf{v}\left(  s\right)  \right)
_{s\in S}$ is a tight frame with frame constant $c$ if and only if the
following two conditions hold:

\begin{enumerate}
\item \label{CorPro.2(1)}$k\left(  t,s\right)  =\overline{k\left(  s,t\right)
}\,$, and

\item \label{CorPro.2(2)}$\displaystyle\sum_{t\in S}k\left(  s_{1},t\right)
k\left(  t,s_{2}\right)  =c^{-1}k\left(  s_{1},s_{2}\right)  $ for all
$s_{1},s_{2}\in S$.
\end{enumerate}
\end{corollary}

\begin{proof}
It is immediate that the stated conditions are equivalent to the matrix $P$ in
(\ref{eqPro.5}) defining a projection when $P=cK$ and $K\left(  t,s\right)
:=\left\langle \,\mathbf{v}\left(  t\right)  \mid\mathbf{v}\left(  s\right)
\,\right\rangle $. As a consequence, we see that the two conditions
(\ref{CorPro.2(1)}) and (\ref{CorPro.2(2)}) are a restatement of
(\ref{eqPre.9}), i.e., the definition of a projection. (We shall consider only
orthogonal projections $P$, i.e., operators $P$ where both of the conditions
in (\ref{eqPre.9}) are assumed.)
\end{proof}

\begin{remark}
\label{RemPro.3}As an application, we are now able to verify the assertion in
Example \textup{\ref{ExaPre.2}}. First note that the Gram matrix of
\textup{(\ref{eqPre.7})} is%
\begin{equation}
k\left(  s_{1},s_{2}\right)  =\cos\left(  \frac{2\pi\left(  s_{1}%
-s_{2}\right)  }{n}\right)  ,\qquad s_{1},s_{2}\in S_{n}. \label{eqPro.9}%
\end{equation}
Hence property \textup{(\ref{CorPro.2(1)})} in the corollary is immediate. To
verify \textup{(\ref{CorPro.2(2)})}, note that
\begin{align*}
\sum_{t\in S_{n}}k\left(  s_{1},t\right)  k\left(  t,s_{2}\right)   &
=\sum_{t\in S_{n}}\cos\left(  \frac{2\pi\left(  s_{1}-t\right)  }{n}\right)
\cos\left(  \frac{2\pi\left(  t-s_{2}\right)  }{n}\right) \\
&  =\frac{1}{2}\sum_{t\in S_{n}}\left(  \cos\left(  \frac{2\pi\left(
s_{1}-s_{2}\right)  }{n}\right)  +\cos\left(  \frac{2\pi\left(  s_{1}%
+s_{2}-2t\right)  }{n}\right)  \right) \\
&  =\frac{n}{2}\sum_{t\in S_{n}}\cos\left(  \frac{2\pi\left(  s_{1}%
-s_{2}\right)  }{n}\right)  =c^{-1}k\left(  s_{1},s_{2}\right)
\end{align*}
with $c=\frac{2}{n}$.

In the last step, we used the trigonometric formula%
\[
\sum_{l=1}^{n}\cos\left(  \frac{4\pi l}{n}-\theta\right)  =0.
\]
\textup{(}We refer to Figures \textup{\ref{FigExaPre.2.1}} to
\textup{\ref{FigExaPre.2.2}} for simple illustrations.\textup{)}\qed
\end{remark}

\begin{corollary}
\label{CorPro.4}Let the pair $\left(  S,\mathcal{H}\right)  $ be as in
Corollary \textup{\ref{CorPro.2}} and Theorem \textup{\ref{ThmPre.4}}. Suppose
$\left\{  \,\mathbf{v}\left(  s\right)  \mid s\in S\,\right\}  $ is a tight
frame in $\mathcal{H}$ with frame constant $c$. Then every $f$ in
$\mathcal{H}$ has the representation%
\begin{equation}
f=c\sum_{s\in S}\left\langle \,\mathbf{v}\left(  s\right)  \mid
f\,\right\rangle \mathbf{v}\left(  s\right)  , \label{eqPro.10}%
\end{equation}
where the sum converges in the norm of the Hilbert space $\mathcal{H}$.
\end{corollary}

\begin{proof}
In view of the argument in the proof of Corollary \ref{CorPro.2}, we may
reduce to the case where the frame constant $c$ is one, i.e., $c=1$. (In
general, the Gram matrix $K\left(  t,s\right)  =\left\langle \,\mathbf{v}%
\left(  t\right)  \mid\mathbf{v}\left(  s\right)  \,\right\rangle
_{\mathcal{H}}$ satisfies $K=c^{-1}P$ where $P$ is a projection in the Hilbert
space $\ell^{2}\left(  S\right)  $.) The reduction to the case $c=1$ means
that%
\begin{equation}
V\colon\mathcal{H}\ni f\longmapsto\left(  \left\langle \,\mathbf{v}\left(
s\right)  \mid f\,\right\rangle \right)  _{s\in S}\in\ell^{2}\left(  S\right)
\label{eqPro.11}%
\end{equation}
is isometric. The projection $P:=VV^{\ast}$ is multiplication by the matrix
\[
\left(  \left\langle \,\mathbf{v}\left(  t\right)  \mid\mathbf{v}\left(
s\right)  \,\right\rangle \right)  _{s,t\in S}.
\]
Since $P$ is the projection onto the range of $V$ in (\ref{eqPro.11}) we
conclude that%
\begin{equation}
\left\langle \,\mathbf{v}\left(  t\right)  \mid f\,\right\rangle =\sum_{s\in
S}\left\langle \,\mathbf{v}\left(  t\right)  \mid\mathbf{v}\left(  s\right)
\,\right\rangle \left\langle \,\mathbf{v}\left(  s\right)  \mid
f\,\right\rangle , \label{eqProprojection}%
\end{equation}
where the convergence is in $\ell^{2}\left(  S\right)  $. But the vectors
$\left\{  \,\mathbf{v}\left(  t\right)  \mid t\in S\,\right\}  $ span a dense
subspace in $\mathcal{H}$, so we get the desired formula%
\begin{equation}
f=\sum_{s\in S}\left\langle \,\mathbf{v}\left(  s\right)  \mid
f\,\right\rangle \mathbf{v}\left(  s\right)  , \label{eqProformula}%
\end{equation}
now referring to the norm and the inner product in $\mathcal{H}$. Recall that
on the range of $V$, the respective inner products of $\mathcal{H}$ and of
$\ell^{2}\left(  S\right)  $ coincide.
\end{proof}

\begin{corollary}
\label{CorPro.5}Let the pair $\left(  S,\mathcal{H}\right)  $ be as in
Corollary \textup{\ref{CorPro.4}}, and suppose that
\[
\left\{  \,\mathbf{v}\left(  s\right)  \mid s\in S\,\right\}
\]
is a tight frame for $\mathcal{H}$ with frame constant $c$. Further suppose
that some $\xi\colon S\rightarrow\mathbb{C}$ satisfies%
\begin{equation}
f=c\sum_{s\in S}\xi_{s}\mathbf{v}\left(  s\right)  , \label{eqPro.12}%
\end{equation}
where the sum is convergent in $\mathcal{H}$. Then%
\begin{equation}
\sum_{s\in S}\left\vert \xi_{s}\right\vert ^{2}\geq\sum_{s\in S}\left\vert
\left\langle \,\mathbf{v}\left(  s\right)  \mid f\,\right\rangle \right\vert
^{2}. \label{eqPro.13}%
\end{equation}

\end{corollary}

\begin{proof}
With the assumptions in the corollary, apply the mapping $V$ from
(\ref{eqPro.11}) to both sides in (\ref{eqPro.12}). We get the formula%
\begin{equation}
\left(  P\xi\right)  _{t}=\left(  \left\langle \,\mathbf{v}\left(  t\right)
\mid f\,\right\rangle _{\mathcal{H}}\right)  \in\ell^{2}\left(  S\right)  .
\label{eqCorPro.5proof(1)}%
\end{equation}
Hence%
\begin{equation}
\xi=P\xi+\left(  I-P\right)  \xi\label{eqCorPro.5proof(2)}%
\end{equation}
is an orthogonal decomposition, i.e.,%
\begin{equation}
\left\Vert \xi\right\Vert _{\ell^{2}}^{2}=\left\Vert P\xi\right\Vert
_{\ell^{2}}^{2}+\left\Vert \left(  I-P\right)  \xi\right\Vert _{\ell^{2}}^{2},
\label{eqCorPro.5proof(3)}%
\end{equation}
and the conclusion (\ref{eqPro.13}) is immediate.
\end{proof}

\section{\label{Sha}Shannon's example}

In the Hilbert space $L^{2}\left(  \mathbb{R}\right)  $ we will consider the
usual Fourier transform%
\begin{equation}
\hat{f}\left(  t\right)  :=\int_{\mathbb{R}}e^{-i2\pi tx}f\left(  x\right)
\,dx, \label{eqSha.1}%
\end{equation}
where convergence is understood in the sense of $L^{2}$. The familiar
interpolation formula of Shannon \cite{Ash90} applies to band-limited
functions, i.e., to functions $f$ on $\mathbb{R}$ such that the Fourier
transform $\hat{f}$ is of compact support. For the present purpose, pick the
following normalization%
\begin{equation}
\operatorname*{supp}\left(  \hat{f}\right)  \subset\left[  -\frac{1}{2}%
,\frac{1}{2}\right]  , \label{eqSha.2}%
\end{equation}
and let $\mathcal{H}$ denote the subspace in $L^{2}\left(  \mathbb{R}\right)
$ defined by this support condition. In particular, $\mathcal{H}$ is the range
of the projection operator $P$ in $L^{2}\left(  \mathbb{R}\right)  $ defined
by%
\begin{equation}
\left(  Pf\right)  \left(  x\right)  :=\int_{-\frac{1}{2}}^{\frac{1}{2}%
}e^{i2\pi xt}\hat{f}\left(  t\right)  \,dt. \label{eqSha.3}%
\end{equation}
Shannon's interpolation formula applies to $f\in\mathcal{H}$, and it reads:%
\begin{equation}
f\left(  x\right)  =\sum_{n\in\mathbb{Z}}f\left(  n\right)  \frac{\sin
\pi\left(  x-n\right)  }{\pi\left(  x-n\right)  }. \label{eqSha.4}%
\end{equation}

\begin{definitions}
\label{DefSha.1}Let $S\subset\mathbb{R}$, and set%
\begin{equation}
\mathbf{v}\left(  s\right)  \left(  x\right)  :=\mathbf{v}\left(  s,x\right)
=\frac{\sin\pi\left(  x-s\right)  }{\pi\left(  x-s\right)  },\qquad s\in S.
\label{eqSha.5}%
\end{equation}
Hence if we take as index set $S:=\mathbb{Z}$, then we may observe that the
functions on the right-hand side in Shannon's formula (\ref{eqSha.4}) are
$\mathbf{v}\left(  n\right)  $ frame vectors, $n\in\mathbb{Z}$. We shall be
interested in other index sets $S$, so-called sets of \emph{sampling points}.
We shall view the functions $\mathbf{v}\left(  s\right)  $ as vectors in
$\mathcal{H}$. The inner product in $\mathcal{H}$ will be that which is
induced from $L^{2}\left(  \mathbb{R}\right)  $, i.e.,%
\begin{equation}
\left\langle \,f_{1}\mid f_{2}\,\right\rangle :=\int_{\mathbb{R}}%
\overline{f_{1}\left(  x\right)  }\,f_{2}\left(  x\right)  \,dx.
\label{eqSha.6}%
\end{equation}

\end{definitions}

The following is well known but is included as an application of Corollary
\ref{CorPro.2}. It is also an example of a pair $\left(  S,\mathcal{H}\right)
$ where $S$ is \emph{infinite}.

\begin{proposition}
\label{ProSha.2}Let $S\subset\mathbb{R}$ be a fixed discrete subgroup, and
assume that $\mathbb{Z}\subset S$. Then $\left\{  \,\mathbf{v}\left(
s\right)  \mid s\in S\,\right\}  $ is a tight frame in $\mathcal{H}$ if and
only if the group index $\left(  S\mathop{:}\mathbb{Z}\right)  $ is finite,
and in that case the frame constant $c$ is $c=\left(  S\mathop{:}\mathbb{Z}%
\right)  ^{-1}$. For the Gram matrix, we have:%
\begin{equation}
K\left(  s_{1},s_{2}\right)  =%
\begin{cases}
\displaystyle\frac{\sin\pi\left(  s_{1}-s_{2}\right)  }{\pi\left(  s_{1}%
-s_{2}\right)  } & \text{for }s_{1},s_{2}\in S,\;s_{1}\neq s_{2},\\
1 & \text{if }s_{1}=s_{2}.
\end{cases}
\label{eqSha.7}%
\end{equation}

\end{proposition}

\begin{proof}
Formula (\ref{eqSha.7}) for the Gram matrix follows from Fourier transform and
the following computation of the inner products:%
\begin{equation}
\left\langle \,\mathbf{v}\left(  s_{1}\right)  \mid\mathbf{v}\left(
s_{2}\right)  \,\right\rangle _{\mathcal{H}}=\int_{\mathbb{R}}\frac{\sin
\pi\left(  x-s_{1}\right)  }{\pi\left(  x-s_{1}\right)  }\,\frac{\sin\pi\left(
x-s_{2}\right)  }{\pi\left(  x-s_{2}\right)  }\,dx. \label{eqProSha.2proof1}%
\end{equation}
A second computation shows that $K\left(  s_{1},s_{2}\right)  =\left\langle
\,\mathbf{v}\left(  s_{1}\right)  \mid\mathbf{v}\left(  s_{2}\right)
\,\right\rangle $ satisfies the two conditions (\ref{CorPro.2(1)}) and
(\ref{CorPro.2(2)}) of Corollary \ref{CorPro.2}, i.e.,%
\begin{equation}
K\left(  s_{1},s_{2}\right)  =K\left(  s_{2},s_{1}\right)  ,
\label{eqProSha.2proof2}%
\end{equation}
and%
\begin{equation}
\sum_{t\in S}K\left(  s_{1},t\right)  K\left(  t,s_{2}\right)  =\left(
S\mathop{:}\mathbb{Z}\right)  K\left(  s_{1},s_{2}\right)  \text{\qquad for
all }s_{1},s_{2}\in S. \label{eqProSha.2proof3}%
\end{equation}

The argument behind this formula uses the known fact that the functions
\[
\left\{  \,\mathbf{v}\left(  n\right)  \mid n\in\mathbb{Z}\,\right\}
\]
in (\ref{eqSha.4}) form an ONB in $\mathcal{H}$; in particular, that%
\begin{equation}
\left\langle \,\mathbf{v}\left(  n_{1}\right)  \mid\mathbf{v}\left(
n_{2}\right)  \,\right\rangle _{\mathcal{H}}=\delta_{n_{1},n_{2}}\text{\qquad
for }n_{1},n_{2}\in\mathbb{Z}. \label{eqProSha.2proof4}%
\end{equation}
Since $\mathbb{Z}\subset S$, we get the following summation formula:%
\begin{align*}
\sum_{t\in S}K\left(  s_{1},t\right)  K\left(  t,s_{2}\right)   &  =\sum_{t\in
S}\left\langle \,\mathbf{v}\left(  s_{1}\right)  \mid\mathbf{v}\left(
t\right)  \,\right\rangle _{\mathcal{H}}\left\langle \,\mathbf{v}\left(
t\right)  \mid\mathbf{v}\left(  s_{2}\right)  \,\right\rangle _{\mathcal{H}}\\
&  =\sum_{k\in S/\mathbb{Z}}\,\sum_{n\in\mathbb{Z}}\left\langle \,\mathbf{v}%
\left(  s_{1}\right)  \mid\mathbf{v}\left(  k+n\right)  \,\right\rangle
_{\mathcal{H}}\left\langle \,\mathbf{v}\left(  k+n\right)  \mid\mathbf{v}%
\left(  s_{2}\right)  \,\right\rangle _{\mathcal{H}}\\
&  =\sum_{k\in S/\mathbb{Z}}\,\sum_{n\in\mathbb{Z}}\frac{\sin\pi\left(
s_{1}-k-n\right)  }{\pi\left(  s_{1}-k-n\right)  }\,\frac{\sin\pi\left(
n-\left(  s_{2}-k\right)  \right)  }{\pi\left(  n-\left(  s_{2}-k\right)
\right)  }\\
&  =\sum_{k\in S/\mathbb{Z}}\frac{\sin\pi\left(  s_{1}-s_{2}\right)  }%
{\pi\left(  s_{1}-s_{2}\right)  }\\
&  =\left(  S\mathop{:}\mathbb{Z}\right)  \left\langle \,\mathbf{v}\left(
s_{1}\right)  \mid\mathbf{v}\left(  s_{2}\right)  \,\right\rangle
_{\mathcal{H}}=\left(  S\mathop{:}\mathbb{Z}\right)  K\left(  s_{1}%
,s_{2}\right)  ,
\end{align*}
which is the desired formula (\ref{eqProSha.2proof3}). The remaining
conclusions in the proposition now follow immediately from Corollary
\ref{CorPro.2}.
\end{proof}

\begin{remark}
\label{RemSha.pound3}The significance of using a larger subgroup $S$, i.e.,
$\mathbb{Z}\subset S$, in a modified version of Shannon's interpolation
formula \textup{(\ref{eqSha.4})} is that a larger \textup{(}%
discrete\/\textup{)} group represents \textquotedblleft\emph{oversampling\/}%
\textquotedblright. However, note that the oversampling changes the frame constant.
\end{remark}

As a contrast showing stability, we now recast a result on oversampling from
\cite{BJMP05} in the present context. It is for tight frames of \emph{wavelet
bases} in $L^{2}\left(  \mathbb{R}\right)  $, and it represents an instance of
stability: a case when oversampling leaves invariant the frame constant.

\begin{proposition}
\label{ProSha.pound4}Let $\psi\in L^{2}\left(  \mathbb{R}\right)  $, and
suppose that the family
\begin{equation}
\psi_{j,k}\left(  x\right)  :=2^{\,j/2}\psi\left(  2^{\,j}x-k\right)  ,\qquad
j,k\in\mathbb{Z}, \label{eqSha.9}%
\end{equation}
is a Parseval frame in $L^{2}\left(  \mathbb{R}\right)  $. Let $p\in
\mathbb{N}$ be odd, $p>1$, and set%
\begin{equation}
\tilde{\psi}_{p}\left(  x\right)  :=\frac{1}{p}\,\psi\left(  \frac{x}{p}\right)
, \label{eqSha.10}%
\end{equation}
and%
\begin{equation}
\tilde{\psi}_{p,j,k}\left(  x\right)  :=2^{\,j/2}\,\tilde{\psi}_{p}\left(
2^{\,j}x-k\right)  ,\qquad j,k\in\mathbb{Z}. \label{eqSha.11}%
\end{equation}
Then the \textquotedblleft oversampled\textquotedblright\ family
\textup{(\ref{eqSha.11})} is again a Parseval frame in the Hilbert space
$L^{2}\left(  \mathbb{R}\right)  $.
\end{proposition}

\begin{proof}
We refer the reader to the argument in Section 2 of \cite{BJMP05}, and to Example \ref{ExaMorNew.6} below.
\end{proof}

\section{\label{Sym}Symmetries}

A basic fact of Hilbert space is that permutations induce unitary operators.
By this we mean that a permutation of the vectors in an orthonormal basis for
a Hilbert space $\mathcal{H}$ induces a unitary transformation $U$ in
$\mathcal{H}$. In this section we explore generalizations of this to frames.

\begin{definition}
\label{DefSym.1}Let $S$ be a set and let $k\colon S\times S\rightarrow
\mathbb{C}$ be a function. We say that $k$ is \emph{positive definite}, or
more precisely positive semidefinite, if the following holds for all finite
sums:%
\begin{equation}
\sum_{s\in S}\sum_{t\in S}\bar{\xi}_{s}\xi_{t}k\left(  s,t\right)  \geq0.
\label{eqSym.1}%
\end{equation}

\end{definition}

While the following result is known, it is not readily available in the
literature, at least not precisely in the form in which we need it. Thus we
include here its statement to save readers from having to track it down.

\begin{theorem}
\label{ThmSym.2}\textup{(Kolmogorov, Parthasarathy--Schmidt \cite{PaSc72})}

\begin{enumerate}
\item \label{ThmSym.2(1)}Let $k\colon S\times S\rightarrow\mathbb{C}$ be
positive definite. Then there are a Hilbert space $\mathcal{H}$ and vectors
$\left\{  \,\mathbf{v}\left(  s\right)  \mid s\in S\,\right\}  \subset
\mathcal{H}$ such that%
\begin{equation}
k\left(  s,t\right)  =\left\langle \,\mathbf{v}\left(  s\right)
\mid\mathbf{v}\left(  t\right)  \,\right\rangle _{\mathcal{H}},\qquad s,t\in
S, \label{eqSym.2}%
\end{equation}
and $\mathcal{H}$ is the closed linear span of $\left\{  \,\mathbf{v}\left(
s\right)  \mid s\in S\,\right\}  $.

\item \label{ThmSym.2(2)}Let $\pi\colon S\rightarrow S$ be a bijection. Then%
\begin{equation}
U_{\pi}\colon\mathbf{v}\left(  s\right)  \longmapsto\mathbf{v}\left(
\pi\left(  s\right)  \right)  \label{eqSym.3}%
\end{equation}
extends to a unitary operator in $\mathcal{H}$ \textup{(}$=:\mathcal{H}\left(
k\right)  $\textup{)} if and only if%
\begin{equation}
k\left(  \pi\left(  s\right)  ,\pi\left(  t\right)  \right)  =k\left(
s,t\right)  \text{\qquad for all }s,t\in S. \label{eqSym.4}%
\end{equation}

\item \label{ThmSym.2(3)}If $\lambda\colon S\rightarrow\mathbb{T}=\left\{
\,z\in\mathbb{C}\mid\left\vert z\right\vert =1\,\right\}  $ and if $\pi$ is as
in \textup{(\ref{ThmSym.2(2)})}, then%
\begin{equation}
U_{\pi,\lambda}\colon\mathbf{v}\left(  s\right)  \longmapsto\lambda\left(
s\right)  \mathbf{v}\left(  \pi\left(  s\right)  \right)  \label{eqSym.5}%
\end{equation}
extends to a unitary operator in $\mathcal{H}$ if and only if%
\begin{equation}
\overline{\lambda\left(  s\right)  }\,\lambda\left(  t\right)  k\left(
\pi\left(  s\right)  ,\pi\left(  t\right)  \right)  =k\left(  s,t\right)
\text{\qquad for all }s,t\in S. \label{eqSym.6}%
\end{equation}

\end{enumerate}
\end{theorem}

\begin{proof}
The reader is referred to \cite{PaSc72} for details.
\end{proof}

In fact, there are important applications of Theorem \ref{ThmSym.2} to frames
with continuous index set. The best known (also included in \cite{PaSc72}) is
the example when $\mathbf{v}\left(  t\right)  $ is a Wiener process, i.e., a
mathematical realization of Einstein's Brownian motion.

\begin{remark}
\label{RemSymNew.pound}As an application of Theorem \textup{\ref{ThmSym.2}(\ref{ThmSym.2(2)})}, note that in each of
the Figures \textup{\ref{FigExaPre.2.1}} and \textup{\ref{FigExaPre.2.2}}, the frames are unitarily equivalent as the angle $\theta$
varies. But they are inequivalent as $n$ varies from $3$ to $4$. More generally,
turning to Example \textup{\ref{ExaPre.2}}, if $n$ is fixed and a new system $\left( \mathbf{v}_{\theta}\left(s\right)\right)$ is
defined by translating the argument in \textup{(\ref{eqPre.7})} by $\theta$, then the assignment
$\mathbf{v}\left(s\right) \to \mathbf{v}_{\theta}\left(s\right)$ extends to a unitary operator in the two-dimensional
Hilbert space $\mathcal{H}$. But as $n$ varies, we get families that are not unitarily
equivalent. In fact, it follows from Remark \textup{\ref{RemPro.3}} that the examples from
Figure \textup{\ref{FigExaPre.2.1}} 
are not equivalent to those in Figure \textup{\ref{FigExaPre.2.2}}, even when equivalence is defined in
the less restrictive sense of \textup{(\ref{eqSym.5})} in Theorem \textup{\ref{ThmSym.2}(\ref{ThmSym.2(3)})}, i.e., allowing a
phase factor in the transformation of the respective systems of frame
vectors.
\end{remark}

For simplicity, return here to the case when the set $S$ is assumed countable
and discrete.

\begin{definition}
\label{DefSym.3}\textup{(Following \cite{Nel57,Nel59}.)}\enspace A closed
subspace $\mathcal{H}$ in $\ell^{2}\left(  S\right)  $ is said to be a
\emph{reproducing kernel subspace} if for every $s\in S$, the mapping%
\begin{equation}
\mathcal{H}\ni f\longmapsto f\left(  s\right)  \in\mathbb{C} \label{eqSym.7}%
\end{equation}
is continuous. Note that by Riesz's lemma this means that there is for each
$s$ a unique element $\mathbf{v}\left(  s\right)  \in\mathcal{H}$ such that%
\begin{equation}
f\left(  s\right)  =\left\langle \,\mathbf{v}\left(  s\right)  \mid
f\,\right\rangle . \label{eqSym.8}%
\end{equation}
And by Schwarz's inequality we get%
\begin{equation}
\left\vert f\left(  s\right)  \right\vert \leq\left\Vert \mathbf{v}\left(
s\right)  \right\Vert _{\mathcal{H}}\left\Vert f\right\Vert _{\mathcal{H}}.
\label{eqSym.9}%
\end{equation}

\end{definition}

\begin{theorem}
\label{ThmSym.4}Let the pair $\left(  S,\mathcal{H}\right)  $ be as in the
statement of Theorem \textup{\ref{ThmPre.4}}. Specifically, we assume that
there is a tight frame $\left(  \mathbf{v}\left(  s\right)  \right)  _{s\in
S}$ for $\mathcal{H}$ with frame constant $c$.

Then it follows that $\mathcal{H}$ is a reproducing kernel Hilbert space, and
that%
\begin{equation}
\left\vert f\left(  s\right)  \right\vert \leq k\left(  s,s\right)
^{1/2}\left\Vert f\right\Vert _{\mathcal{H}}\text{\qquad for all }s\in S\text{
and }f\in\mathcal{H}, \label{eqSym.10}%
\end{equation}
where $k\left(  s,t\right)  =\left\langle \,\mathbf{v}\left(  s\right)
\mid\mathbf{v}\left(  t\right)  \,\right\rangle $ is the Gram matrix of
$\mathcal{H}$.
\end{theorem}

\begin{proof}
Let $\left(  S,\mathcal{H}\right)  $ be as stated. Then by Corollary
\ref{CorPro.4}, we have the representation%
\begin{equation}
f=c\sum_{s\in S}\left\langle \,\mathbf{v}\left(  s\right)  \mid
f\,\right\rangle \mathbf{v}\left(  s\right)  , \label{eqSym.11}%
\end{equation}
referring to the isometric embedding $\mathcal{H}\underset{\simeq
\rule[-0.5ex]{0pt}{0.5ex}}{\hooklongrightarrow}\ell^{2}\left(  S\right)  $. The Gram matrix of $\left(
\mathbf{v}\left(  s\right)  \right)  _{s\in S}$ induces an operator $K$, and
$P=cK$ is the projection of $\ell^{2}\left(  S\right)  $ onto the subspace
$\mathcal{H}$. Moreover, $\mathbf{v}\left(  s\right)  =P\left(  \delta
_{s}\right)  $ for all $s\in S$. Now apply both sides of (\ref{eqSym.11}) to
some point $t\in S$. Via the isometric embedding, we know that%
\begin{equation}
\mathcal{H}=\left\{  \,f\in\ell^{2}\left(  S\right)  \mid Pf=f\,\right\}  .
\label{eqThmSym.4proof1}%
\end{equation}
Hence, if $f\in\mathcal{H}$, we get $f\left(  t\right)  =\left\langle
\,\delta_{t}\mid f\,\right\rangle =\left\langle \,\delta_{t}\mid
Pf\,\right\rangle =\left\langle \,P\delta_{t}\mid f\,\right\rangle
=\left\langle \,\mathbf{v}\left(  t\right)  \mid f\,\right\rangle $; and by
(\ref{eqSym.11}),%
\begin{equation}
f\left(  t\right)  =c\sum_{s\in S}\left\langle \,\mathbf{v}\left(  s\right)
\mid f\,\right\rangle \,\underbrace{\left\langle \,\mathbf{v}\left(  t\right)
\mid\mathbf{v}\left(  s\right)  \,\right\rangle }_{k\left(  t,s\right)  }.
\label{eqThmSym.4proof2}%
\end{equation}
An application of Corollary \ref{CorPro.2}, and of Schwarz for $\ell
^{2}\left(  S\right)  $, now yields%
\begin{align*}
\left\vert f\left(  t\right)  \right\vert  &  \leq c\left(  \sum_{s\in
S}\left\vert \left\langle \,\mathbf{v}\left(  s\right)  \mid f\,\right\rangle
\right\vert ^{2}\right)  ^{1/2}\left(  \sum_{s\in S}\left\vert k\left(
t,s\right)  \right\vert ^{2}\right)  ^{1/2}\\
&  =c\cdot c^{-1/2}\left\Vert f\right\Vert _{\mathcal{H}}c^{-1/2}\,k\left(
t,t\right)  ^{1/2}\\
&  =\left\Vert f\right\Vert _{\mathcal{H}}k\left(  t,t\right)  ^{1/2}.
\settowidth{\qedskip}{$\displaystyle\left\vert f\left(  t\right)  \right\vert
=\left\Vert f\right\Vert _{\mathcal{H}}k\left(  t,t\right)  ^{1/2}.$}\settowidth{\qedadjust}{$\displaystyle\left\vert f\left(  t\right)  \right\vert  \leq c\left(  \sum_{s\in
S}\left\vert \left\langle \,\mathbf{v}\left(  s\right)  \mid f\,\right\rangle
\right\vert ^{2}\right)  ^{1/2}\left(  \sum_{s\in S}\left\vert k\left(
t,s\right)  \right\vert ^{2}\right)  ^{1/2}$}\addtolength{\qedskip}{-0.5\qedadjust}\addtolength{\qedskip}{-0.5\textwidth}\rlap{\makebox[-\qedskip]{\qed}}
\end{align*}
\renewcommand{\qed}{\relax}
\end{proof}

\begin{remark}
\label{RemSym.5}It is clear that the converse to the theorem also holds: If
$\mathcal{H}\underset{\simeq\rule[-0.5ex]{0pt}{0.5ex}}{\hooklongrightarrow}\ell^{2}\left(  S\right)  $
is a reproducing kernel Hilbert space, then for each $s\in S$, and
$f\in\mathcal{H}$, $f\left(  s\right)  =\left\langle \,\mathbf{v}\left(
s\right)  \mid f\,\right\rangle _{\mathcal{H}}$ holds for a unique family
$\left(  \mathbf{v}\left(  s\right)  \right)  _{s\in S}$ in $\mathcal{H}$.
This family will be a tight frame with frame constant one. If $c\in\left(
0,1\right)  $, then $\mathbf{w}_{c}\left(  s\right)  :=c^{-1/2}\mathbf{v}%
\left(  s\right)  $ will define a tight frame with frame constant $c$.
\end{remark}

\section{\label{Mor}More general frames}

Since the condition (\ref{eqPre.6}) which defines tight frames is rather
rigid, it is of interest to consider how it can be relaxed in a way which
still makes it useful.

As before, we will consider a pair $\left(  S,\mathcal{H}\right)  $, where $S$
is a fixed countable set, and where $\mathcal{H}$ is a Hilbert space

\begin{definitions}
\label{DefMor.1}A system of vectors $\left(  \mathbf{v}\left(  s\right)
\right)  _{s\in S}$ in $\mathcal{H}$ is called a \emph{frame} for
$\mathcal{H}$ if there are constants $0<A_{1}\leq A_{2}<\infty$ such that%
\begin{equation}
A_{1}\left\Vert f\right\Vert ^{2}\leq\sum_{s\in S}\left\vert \left\langle
\,\mathbf{v}\left(  s\right)  \mid f\,\right\rangle \right\vert ^{2}\leq
A_{2}\left\Vert f\right\Vert ^{2}\text{\qquad for all }f\in\mathcal{H}.
\label{eqMor.1}%
\end{equation}

It is called a \emph{Bessel sequence} if only the estimate on the right-hand
side in \textup{(\ref{eqMor.1})} is assumed, i.e., if for some finite constant
$A$,%
\begin{equation}
\sum_{s\in S}\left\vert \left\langle \,\mathbf{v}\left(  s\right)  \mid
f\,\right\rangle \right\vert ^{2}\leq A\left\Vert f\right\Vert ^{2}%
\text{\qquad for all }f\in\mathcal{H}. \label{eqMor.2}%
\end{equation}

Recall that when \textup{(\ref{eqMor.2})} is assumed, then the analysis
operator $V=V_{\left(  \mathbf{v}\left(  s\right)  \right)  }$ given by%
\begin{equation}
\mathcal{H}\ni f\overset{V}{\longmapsto}\left(  \left\langle \,\mathbf{v}%
\left(  s\right)  \mid f\,\right\rangle \right)  _{s\in S}\in\ell^{2}\left(
S\right)  \label{eqMor.3}%
\end{equation}
is well defined and bounded. Hence, the adjoint operator $V^{\ast}\colon
\ell^{2}\left(  S\right)  \rightarrow\mathcal{H}$ is bounded as well, and
\begin{equation}
V^{\ast}\left(  \xi_{s}\right)  =\sum_{s\in S}\xi_{s}\,\mathbf{v}\left(
s\right)  \text{\qquad for all }\left(  \xi_{s}\right)  \in\ell^{2}\left(
S\right)  , \label{eqMor.4}%
\end{equation}
where the sum on the right-hand side in \textup{(\ref{eqMor.4})} is convergent
in $\mathcal{H}$ for all $\left(  \xi_{s}\right)  \in\ell^{2}\left(  S\right)
$.

The tight frames for which $A_{1}=A_{2}=1$ are called \emph{Parseval frames}.
\end{definitions}

\begin{theorem}
\label{ThmMor.2}Let $\left(  S,\mathcal{H}\right)  $ be as above, and let
$\left(  \mathbf{v}\left(  s\right)  \right)  _{s\in S}$ be a Bessel sequence
with Bessel constant $A$.

\begin{enumerate}
\item \label{ThmMor.2(1)}Then the closed span $\mathcal{H}_{\operatorname*{in}%
}$ of $\left(  \mathbf{v}\left(  s\right)  \right)  _{s\in S}$ contains a
derived Parseval frame.

\item \label{ThmMor.2(2)}The derived Parseval frame is a Parseval frame for
$\mathcal{H}$ if and only if $\left(  \mathbf{v}\left(  s\right)  \right)
_{s\in S}$ is a frame for $\mathcal{H}$, i.e., iff $\mathcal{H}%
_{\operatorname*{in}}=\mathcal{H}$.

\item \label{ThmMor.2(3)}In the general case when $\left(  \mathbf{v}\left(
s\right)  \right)  _{s\in S}$ is a Bessel sequence, the operator $W:=V\left(
V^{\ast}V\right)  ^{-1/2}$ is well defined and isometric on $\mathcal{H}%
_{\operatorname*{in}}$, and%
\begin{equation}
\mathbf{w}\left(  s\right)  :=\left(  V^{\ast}V\right)  ^{-1/2}\mathbf{v}%
\left(  s\right)  ,\qquad s\in S, \label{eqMor.5}%
\end{equation}
is a Parseval frame in $\mathcal{H}_{\operatorname*{in}}$.
\end{enumerate}
\end{theorem}

\begin{proof}
It is easy to see that both $V$ and $V^{\ast}$ are bounded; $V$ is everywhere
defined on $\mathcal{H}$, and $V^{\ast}$ everywhere defined on $\ell
^{2}\left(  S\right)  $. For the operator norms, we have $\left\Vert
V\right\Vert =\left\Vert V^{\ast}\right\Vert =\left\Vert V^{\ast}V\right\Vert
^{1/2}\leq\sqrt{A}$. The \emph{initial space} of $V$, $\mathcal{H}\left(
V\right)  $, is defined as%
\begin{equation}
\mathcal{H}\ominus\left\{  \,f\in\mathcal{H}\mid Vf=0\,\right\}
=\overline{R\left(  V^{\ast}\right)  }, \label{eqMor.6}%
\end{equation}
where the over-bar stands for norm-closure in $\mathcal{H}$, and where
$R\left(  V^{\ast}\right)  $ denotes the range of $V^{\ast}$.

We let $\mathcal{H}_{\operatorname*{in}}$ denote the closed span of the given
Bessel sequence $\left(  \mathbf{v}\left(  s\right)  \right)  _{s\in S}$. Our
claim is that%
\begin{equation}
\mathcal{H}_{\operatorname*{in}}=\mathcal{H}\left(  V\right)  .
\label{eqMor.7}%
\end{equation}
To see this, note that (\ref{eqMor.4}) implies the inclusion ($\supseteq$) in
(\ref{eqMor.7}). But all finite linear combinations $\sum_{s}\xi_{s}\,%
\mathbf{v}\left(  s\right)  $ are contained in $R\left(  V^{\ast}\right)  $,
and by closure, we get $\mathcal{H}_{\operatorname*{in}}\subseteq
\overline{R\left(  V^{\ast}\right)  }$, which is the second inclusion in
(\ref{eqMor.7}). Hence (\ref{eqMor.7}) holds.

We now apply the polar decomposition (from operator theory), see
\cite{KaRi97}, to the operator $V$. The conclusion is that there is a
\emph{partial isometric} $W$ such that%
\begin{equation}
V=W\left(  V^{\ast}V\right)  ^{1/2}=\left(  VV^{\ast}\right)  ^{1/2}W
\label{eqMor.8}%
\end{equation}
and such that the \emph{initial space} of $W$ is $\mathcal{H}_{\operatorname*{in}%
}$ and the \emph{final space} of $W$ is $\overline{R\left(  V\right)  }$.
Implied in this are the following assertions:

\begin{enumerate}
\renewcommand{\theenumi}{\roman{enumi}}

\item \label{ThmMor.2proof(1)}Both of the operators $W^{\ast}W$ and $WW^{\ast
}$ are projections;

\item \label{ThmMor.2proof(2)}The range of $W^{\ast}W$ is $\mathcal{H}%
_{\operatorname*{in}}$;

\item \label{ThmMor.2proof(3)}The range of $WW^{\ast}$ is $\overline{R\left(
V\right)  }$;

\item \label{ThmMor.2proof(4)}The operator $V^{\ast}V$ is selfadjoint, and
$\left(  V^{\ast}V\right)  ^{1/2}$ is defined by the spectral theorem applied
to $V^{\ast}V$.
\end{enumerate}

Now apply Theorem \ref{ThmSym.4} to the restriction operator%
\begin{equation}
W\colon\mathcal{H}_{\operatorname*{in}}\longrightarrow\ell^{2}\left(
S\right)  . \label{eqMor.9}%
\end{equation}
In view of (\ref{ThmMor.2proof(2)}) above, this is an isometry. Recall
$W^{\ast}W\mathcal{H}=\mathcal{H}_{\operatorname*{in}}$ by (\ref{eqMor.7}) and
the polar decomposition (\ref{eqMor.8}). As a result, we get that the vectors
$\mathbf{w}\left(  s\right)  =W^{\ast}\delta_{s}$, for $s\in S$, form a
Parseval frame for $\mathcal{H}_{\operatorname*{in}}$. But we also have the
formula $V^{\ast}\delta_{s}=\mathbf{v}\left(  s\right)  $ as an application of
(\ref{eqMor.4}). Using now the two facts (\ref{eqMor.9}) and (\ref{eqMor.8}),
we get $\mathbf{w}\left(  s\right)  =W^{\ast}\delta_{s}=\left(  V^{\ast
}V\right)  ^{-1/2}V^{\ast}\delta_{s}=\left(  V^{\ast}V\right)  ^{-1/2}%
\mathbf{v}\left(  s\right)  $, which is the desired formula (\ref{eqMor.5})
from the conclusion in the theorem. The remaining conclusions stated in the
theorem are already implied by the reasoning above.
\end{proof}

\begin{remark}
\label{RemMor.pound3}\textbf{The operators }$V^{\ast}V$\textbf{ and }%
$VV^{\ast}$\textbf{.}\enspace Except for the point $\lambda=0$, the two
operators $V^{\ast}V$ and $VV^{\ast}$ have the same spectrum. But in the
general case, this spectrum could be discrete or continuous; or it could even
be singular.

The fact that the spectrum minus $\left\{  0\right\}  $ is the same follows
from \textup{(\ref{eqMor.8})}. In the \emph{frame case}, it follows from
\textup{(\ref{eqMor.1})} that the spectrum is contained in the interval
$\left[  A_{1},A_{2}\right]  $. But if $\left(  \mathbf{v}\left(  s\right)
\right)  _{s\in S}$ is only a Bessel sequence, \textup{(\ref{eqMor.2})}, the
best that can be said is that the spectrum is contained in $\left[
0,A\right]  $. In fact, the lower bound in $\operatorname*{spectrum}\left(
V^{\ast}V\right)  $ is also the best lower frame bound, referring to
\textup{(\ref{eqMor.1})}.

Computationally, however, the two operators $V^{\ast}V$ and $VV^{\ast}$ are
quite different. The first one operates in $\mathcal{H}$, while the second one
maps $\ell^{2}\left(  S\right)  $ into itself.

The formula
\begin{equation}
V^{\ast}Vf=\sum_{s\in S}\left\langle \,\mathbf{v}\left(  s\right)  \mid
f\,\right\rangle \mathbf{v}\left(  s\right)  \text{\qquad for }f\in\mathcal{H}
\label{eqRemMor.pound3.1}%
\end{equation}
shows that $V^{\ast}V$ serves to decompose vectors in $\mathcal{H}$ relative
to the given frame $\left(  \mathbf{v}\left(  s\right)  \right)  _{s\in S}$.
In contrast, $VV^{\ast}$ has an explicit matrix representation:
\end{remark}

\begin{proposition}
\label{ProMor.pound4}Relative to the ONB $\left(  \delta_{s}\right)  _{s\in S}%
$, $VV^{\ast}$ is simply multiplication with the Gram matrix 
\begin{equation}
\left(
\left\langle \,\mathbf{v}\left(  s\right)  \mid\mathbf{v}\left(  t\right)
\,\right\rangle _{\mathcal{H}}\right)  _{s,t\in S}.
\label{eqMor.11}%
\end{equation}
\end{proposition}

\begin{proof}
To compute the $\left(  s,t\right)  $-matrix entry for the matrix which
represents the operator $VV^{\ast}$, we use the inner product in $\ell
^{2}\left(  S\right)  $ as follows: The $\left(  s,t\right)  $-matrix entry is%
\begin{equation}
\left\langle \,\delta_{s}\mid VV^{\ast}\delta_{t}\,\right\rangle _{\ell^{2}%
}=\left\langle \,V^{\ast}\delta_{s}\mid V^{\ast}\delta_{t}\,\right\rangle
_{\mathcal{H}}=\left\langle \,\mathbf{v}\left(  s\right)  \mid\mathbf{v}%
\left(  t\right)  \,\right\rangle _{\mathcal{H}},
\label{eqProMor.pound4proof1}%
\end{equation}
which is the Gram matrix.
\end{proof}

\begin{example}
\label{ExaMor.5}Example \textup{\ref{ExaPre.2}} revisited for $n=3$.\enspace

In the case of Example \textup{\ref{ExaPre.2}} \textup{(}see also
Fig.\ \textup{\ref{FigExaPre.2.1})} for $n=3$, $VV^{\ast}$ is the $3\times3$
matrix
\[%
\begin{pmatrix}
1 & -\frac{1}{2} & -\frac{1}{2}\\
-\frac{1}{2} & 1 & -\frac{1^{\mathstrut}}{2_{\mathstrut}}\\
-\frac{1}{2} & -\frac{1}{2} & 1
\end{pmatrix}
,
\]
while $V^{\ast}V$ may be represented as
\[%
\begin{pmatrix}
\,\frac{3}{2_{\mathstrut}} & 0\,\\
\,0 & \frac{3}{2}\,
\end{pmatrix}
,
\]
i.e., a $2\times2$ matrix. Moreover, $\operatorname*{spec}\left(  VV^{\ast
}\right)  =\left\{  0,\frac{3}{2}\right\}  $, and $\operatorname*{spec}\left(
V^{\ast}V\right)  =\left\{  \frac{3}{2}\right\}  $.
\end{example}

\begin{corollary}
\label{CorMor.pound6}Consider a system $\left(  \mathcal{H},S,\left(
\mathbf{v}\left(  s\right)  \right)  _{s\in S}\right)  $ as above, and let
\[
V=V_{\left(  \mathbf{v}\left(  s\right)  \right)  }\colon\mathcal{H}%
\longrightarrow\ell^{2}\left(  S\right)
\]
be the associated analysis operator. Then the respective lower and upper frame
estimates \textup{(\ref{eqMor.1})} are equivalent to
\begin{equation}
\operatorname*{spec}\left(  \left(  \left\langle \,\mathbf{v}\left(  s\right)
\mid\mathbf{v}\left(  t\right)  \,\right\rangle \right)  _{\mathstrut
}^{\mathstrut}\right)  \setminus\left\{  0\right\}  \subset\left[  A_{1}%
,A_{2}\right]  , \label{eqMor.12}%
\end{equation}
where $\left(  \left\langle \,\mathbf{v}\left(  s\right)  \mid\mathbf{v}%
\left(  t\right)  \,\right\rangle \right)  _{s,t\in S}$ is the Gram matrix in
\textup{(\ref{eqMor.11})}.
\end{corollary}

\begin{proof}
Since we noted in Remark \ref{RemMor.pound3} that
\begin{equation}
\operatorname*{spec}\left(  V^{\ast}V\right)  \setminus\left\{  0\right\}
=\operatorname*{spec}\left(  VV^{\ast}\right)  \setminus\left\{  0\right\}  ,
\label{eqMor.13}%
\end{equation}
the assertion in (\ref{eqMor.12}) follows immediately from Proposition
\ref{ProMor.pound4}. Specifically, the two estimates in (\ref{eqMor.1}) are
equivalent to the following system of operator inequalities:%
\begin{equation}
A_{1}I_{\mathcal{H}}\leqq V^{\ast}V\leqq A_{2}I_{\mathcal{H}},
\label{eqMor.14}%
\end{equation}
where the ordering \textquotedblleft$\leqq$\textquotedblright\ in
(\ref{eqMor.14}) refers to the familiar ordering on the set of selfadjoint
operators \cite{KaRi97}. Now an application of the spectral theorem
\cite{KaRi97} to $V^{\ast}V$ shows that (\ref{eqMor.14}) is also equivalent to
the set-containment $\operatorname*{spec}\left(  V^{\ast}V\right)
\subseteq\left[  A_{1},A_{2}\right]  $. In fact, the lower estimate in
(\ref{eqMor.1}) holds iff $\operatorname*{spec}\left(  V^{\ast}V\right)
\subseteqq\left[  A_{1},\infty\right)  $, while the upper estimate is
equivalent to the containment $\operatorname*{spec}\left(  V^{\ast}V\right)
\subseteqq\left[  0,A_{2}\right]  $.
\end{proof}

While Corollary \ref{CorMor.pound6} is rather abstract, it has a variety of
specific uses which serve to make it clear how certain initial frame systems
are transformed into more detailed frames having additional structure. A case
in point is the transformation of frames used in the analysis of signals into
frequency bands, for example mutually non-interfering low-pass and high-pass
bands. A convenient tool for accomplishing this is the discrete wavelet
transform. As noted in, for example, \cite{Chr03} and \cite{Jor06a}, the
discrete wavelet transform serves to create nested families of resolution
subspaces (details below) in the Hilbert spaces of sequences which are used in
our representation of time-series, or of speech signals. When more sub-bands
are introduced, this same idea works for the analysis of images.

By a discrete wavelet transform with perfect reconstruction we shall mean
(following \cite{Jor06a}) a Hilbert space $\mathcal{H}$ and a system of
operators $F_{0}$, $F_{1}$ in $\mathcal{H}$ such that
\begin{equation}
F_{i} F_{j}^{*} = \delta_{i,j} I , \quad\sum_{i=0}^{1} F_{i}^{*} F_{i} =
I,\text{\quad and\quad} F_{0} ^{n} f \underset{n \rightarrow\infty
}{\longrightarrow} 0 \text{\qquad for all }f \in\mathcal{H}%
.\label{eqMorNew.star}%
\end{equation}

As further noted in \cite{Jor06a}, each system of quadrature mirror filters
from standard signal processing defines an operator system $\left(
F_{i}\right) $ as in (\ref{eqMorNew.star}).

\begin{remarks}
\label{RemMorNew.pound}%
\setcounter{enump}{0}\renewcommand{\theenump}{\alph{enump}}\refstepcounter
{enump}\label{RemMorNew.pound1}\textup{(\ref{RemMorNew.pound1})}~The subscript
convention for the two operators $F_{0}$ and $F_{1}$ in
\textup{(\ref{eqMorNew.star})} comes from engineering. The index value $i = 0$
corresponds to a chosen low-pass filter followed by downsampling, while
$F_{1}$ is the operation of a high-pass filtering followed by downsampling.
Hence ``low'' is indexed by zero. This index convention is not related to that
of the two frame bounds $A_{i}$ in \textup{(\ref{eqMor.1})}. These two numbers
are usually called $A_{1}$ and $A_{2}$.

\refstepcounter{enump}\label{RemMorNew.pound2}\textup{(\ref{RemMorNew.pound2}%
)}~[Referring to the sum-formula in \textup{(\ref{eqMorNew.star})}.] In
stating our quadrature conditions \textup{(\ref{eqMorNew.star})} we have for
reasons of clarity restricted attention to the simplest case: dyadic and
orthogonal filters. \textup{(}The dyadic case refers to the dyadic wavelets
from Proposition \textup{\ref{ProSha.pound4}} above.\textup{)} But our present
discussion \textup{(}both in Sections \textup{\ref{Sha}} and \textup{\ref{Mor}%
)} easily generalizes to systems with more than two bands, and even to less
restrictive quadrature conditions. In case of more, say $N$ bands, the
indexing of the operators $F_{i}$ is $i = 0, 1, \dots, N -1$. And instead of
using the operator $F_{i}^{*}$ in the sum in \textup{(\ref{eqMorNew.star})} at
the $i$'th place, we may for some applications rather use a second operator
$G_{i}$ \textup{(}not equal to $F_{i}^{*}$\textup{)}; see for example
\cite{JoKr03} for details on these more general systems.
\end{remarks}

\begin{corollary}
\label{CorMorNew.7}Let $\left( \mathbf{v}\left( s\right) \right) _{s \in S}$
be a frame with frame bounds $A_{i}$, $i = 1,2$, in a Hilbert space
$\mathcal{H}$, and let $\left( F_{i}\right) $ be a discrete wavelet transform
system \textup{(}as in \textup{(\ref{eqMorNew.star}))}. Set
\begin{equation}
\mathbf{v}\left( k, s\right)  : = F_{0}^{*\,k}\, F_{1}^{*} \, \mathbf{v}\left(
s\right) ,\qquad s \in S,\; k = 0, 1, 2, \dots.\label{eqMorNew.starstar}%
\end{equation}
Then the subdivided vector system $\left( \mathbf{v}\left( k, s\right) \right)
$ is a frame in $\mathcal{H}$ with the same frame bounds $A_{i}$ as $\left(
\mathbf{v}\left( s\right) \right) $. Moreover the Gramian of $\left(
\mathbf{v}\left( k, s\right) \right) $ has the form $I \otimes G$ where $I$ is
the \textup{(}infinite\/\textup{)} identity matrix, and $G$ is the Gramian for
$\left( \mathbf{v}\left( s\right) \right) $.

     Moreover, if two frames $\left( \mathbf{v}\left( s\right) \right)$ and $\left( \mathbf{v}^{\prime}\left( s\right) \right)$ are unitarily equivalent in
the sense of Theorem \textup{\ref{ThmSym.2}}, then two refined systems $\left( \mathbf{v}\left( k, s\right) \right)$ and $\left( \mathbf{v}^{\prime}\left( k, s\right) \right)$ resulting from \textup{(\ref{eqMorNew.starstar})} are also unitarily equivalent.
\end{corollary}

\begin{proof}
Since in general, by Corollary \ref{CorMor.pound6}, the two frame bounds
(upper and lower) coincide with the spectral bounds for the corresponding
Gramian, we only need to study the Gramian of the new system $\left(
\mathbf{v}\left( k, s\right) \right) $ in (\ref{eqMorNew.starstar}). Recall
that the Gramian for the new system is the matrix of inner products, with the
only modification to the formula in Theorem \ref{ThmSym.4} being that the
row and column indices are now double indices.

The significance of the third condition from (\ref{eqMorNew.star}) is related
to the kind of subspace structure in $\mathcal{H}$ which models
\emph{resolutions} of signals. Since $F_{0} F_{0}^{*} = I $, we get a nested
system of projections
\begin{equation}
\label{eqCorMorNew.7proof1}P_{k}:=F_{0}^{*\,k} F_{0}^{k} , \qquad k=0,1,\dots,
\end{equation}
and
\begin{equation}
\label{eqCorMorNew.7proof2}\cdots P_{k+1}\leq P_{k}\leq\dots\leq P_{1}\leq
P_{0}=I.
\end{equation}
Recall that the ordering of projections coincides with the associated ordering
of the range subspaces $\mathcal{H}_{k}:=P_{k}\mathcal{H}$ in $\mathcal{H}$,
i.e.,
\begin{equation}
\label{eqCorMorNew.7proof3}\cdots\mathcal{H}_{k+1}\subseteq\mathcal{H}%
_{k}\subseteq\dots\subseteq\mathcal{H}_{1}\subseteq\mathcal{H}.
\end{equation}
The third condition in (\ref{eqMorNew.star}) is equivalent to the assertion
that
\begin{equation}
\label{eqCorMorNew.7proof4}\bigcap_{k\in\left\{ 0,1,\dots\right\} }%
\mathcal{H}_{k}=\left\{ \mathbf{0}\right\} .
\end{equation}
Specifically,
\begin{equation}
\label{eqCorMorNew.7proof5}\left\langle \, \mathbf{v}\left( j, s\right)
\mid\mathbf{v}\left( k, t\right) \,\right\rangle = \left\langle \,
F_{0}^{*\,j}\, F_{1}^{*} \, \mathbf{v}\left( s\right) \bigm| F_{0}^{*\,k}\,
F_{1}^{*} \, \mathbf{v}\left( t\right) \,\right\rangle = \delta_{j,k}
\left\langle \, \mathbf{v}\left( s\right) \mid\mathbf{v}\left( t\right)
\,\right\rangle .
\end{equation}
Since the last expression is the tensor product of $I$ with $G$, the result follows.

     The last part of the conclusion in the corollary about preservation of
unitary equivalence in passing to the refinement $\left( \mathbf{v}\left( s\right) \right)\to\left( \mathbf{v}\left( k, s\right) \right)$
follows from this and Theorem \ref{ThmSym.2}.
\end{proof}

\begin{example}
\label{ExaMorNew.6}\textup{(Following the terminology from
(\ref{eqSha.9}) in Proposition \ref{ProSha.pound4}.)} Suppose
some dyadic wavelet $\left( \psi_{j,k}\right) $ comes from a scaling function
$\varphi$. Set $S:= \mathbb{Z}$, and
\begin{equation}
\mathbf{v}\left( k\right) := \varphi\left( x - k\right) , \qquad k
\in\mathbb{Z}.\label{eqExaMorNew.6(1)}%
\end{equation}
Then there are known conditions on such a scaling function $\varphi$ for the
frame estimates \textup{(\ref{eqMor.1})} to hold; see, e.g., \cite{CoDa93}. To
apply the operator system \textup{(\ref{eqMorNew.star})} from above, choose
the Hilbert space $\mathcal{H}$ be the closed subspace in $L^{2}\left(
\mathbb{R}\right) $ spanned by the integral translates $\left\{ \,
\varphi\left(  x - k\right)  \mid k \in\mathbb{Z}\,\right\} $. Following
\cite{Jor06a} and \textup{(\ref{eqMorNew.starstar})} above, we may then
construct operators $F_{0}$, $F_{1}$ and introduce the associated resolution
frame $\left(  \mathbf{v}\left( j, k\right) \right) $, double-indexed as in
Corollary \textup{\ref{CorMorNew.7}} as follows:

Set
\begin{equation}
\mathbf{v}\left( j, k\right)  := F_{0}^{*\,j}\, F_{1}^{*} \, \mathbf{v}\left(
k\right)  = \psi_{ -j - 1, k} \text{\qquad for }j= 0, 1, 2, \dots, \text{ and
}k \in\mathbb{Z}.\label{eqExaMorNew.6(2)}%
\end{equation}
Using the corollary, we then conclude that for each $j$, the functions
$\left\{ \, \mathbf{v}\left( j, k\right)  \mid k \in\mathbb{Z}\,\right\} $
generate the relative complement subspace $\mathcal{H}_{j} \ominus
\mathcal{H}_{j + 1}$ from the nested resolution system we introduced there. By
this we mean that for each $j$, $\mathcal{H}_{j} \ominus\mathcal{H}_{j + 1}$
is the closed linear span of $\left\{ \, \mathbf{v}\left( j, k\right)  \mid k
\in\mathbb{Z}\,\right\} $.
\end{example}


\begin{thebibliography}{CFTW06}

\bibitem[AkGl93]{AkGl93}
N.I. Akhiezer and I.M. Glazman, {\em Theory of {L}inear {O}perators in
  {H}ilbert {S}pace}, Dover, New York, 1993, reprint of the 1961--63
  publication (Frederick Ungar Publishing Co., New York) of a translation by
  Merlynd Nestell; a corrected and augmented Russian edition (``Vishcha
  Shkola'', Kharkov, 1977--78) has also been translated by E.R. Dawson and
  published in the series Monographs and Studies in Mathematics, vol.~9, Pitman
  (Advanced Publishing Program), Boston--London, 1981.

\bibitem[ALTW04]{ALTW04}
A.~Aldroubi, D.~Larson, W.-S. Tang, and E.~Weber, {\em Geometric aspects of
  frame representations of abelian groups}, Trans. Amer. Math. Soc. {\bf 356}
  (2004), 4767--4786.

\bibitem[Aro50]{Aro50}
N.~Aronszajn, {\em Theory of reproducing kernels}, Trans. Amer. Math. Soc. {\bf
  68} (1950), 337--404.

\bibitem[Ash90]{Ash90}
R.B. Ash, {\em Information {T}heory}, Dover, New York, 1990, corrected reprint
  of the original 1965 Interscience/Wiley edition.

\bibitem[BJMP05]{BJMP05}
L.W. Baggett, P.E.T. Jorgensen, K.D. Merrill, and J.A. Packer, {\em
  Construction of {P}arseval wavelets from redundant filter systems}, J. Math.
  Phys. {\bf 46} (2005), no.~8, 083502, 28 pp., doi:10.1063/1.1982768.

\bibitem[BJMP06]{BJMP06}
\bysame
, {\em A non-{MRA}
  {$C^r$} frame wavelet with rapid decay}, Acta Appl. Math. {\bf Online First}
  (2006), 20 pp., doi:10.1007/s10440-005-9011-4.

\bibitem[BeFi03]{BeFi03}
J.J. Benedetto and M.~Fickus, {\em Finite normalized tight frames}, Adv.
  Comput. Math. {\bf 18} (2003), 357--385.

\bibitem[BDP05]{BDP05}
S.~Bildea, D.E. Dutkay, and G.~Picioroaga, {\em M{RA} super-wavelets}, New York
  J. Math. {\bf 11} (2005), 1--19.

\bibitem[CaCh03]{CaCh03}
P.G. Casazza and O.~Christensen, {\em Gabor frames over irregular lattices},
  Adv. Comput. Math. {\bf 18} (2003), 329--344.

\bibitem[CCLV05]{CCLV05}
P.G. Casazza, O.~Christensen, A.M. Lindner, and R.~Vershynin, {\em Frames and
  the {F}eichtinger conjecture}, Proc. Amer. Math. Soc. {\bf 133} (2005),
  1025--1033.

\bibitem[CFTW06]{CFTW06}
P.G. Casazza, M.~Fickus, J.C. Tremain, and E.~Weber, {\em The
  {K}adison--{S}inger problem in mathematics and engineering: a detailed
  account}, Operator Theory, Operator Algebras, and
  Applications (GPOTS 2005) (D.~Han, P.E.T. Jorgensen, and D.~Larson,
  eds.), Contemp. Math., American Mathematical Society, Providence, to
  appear.

\bibitem[CKL04]{CKL04}
P.G. Casazza, G.~Kutyniok, and M.C. Lammers, {\em Duality principles in frame
  theory}, J. Fourier Anal. Appl. {\bf 10} (2004), 383--408.

\bibitem[Chr99]{Chr99}
O.~Christensen, {\em Operators with closed range, pseudo-inverses, and
  perturbation of frames for a subspace}, Canad. Math. Bull. {\bf 42} (1999),
  37--45.

\bibitem[Chr03]{Chr03}
\bysame
, {\em An {I}ntroduction to {F}rames and {R}iesz {B}ases},
  Applied and Numerical Harmonic Analysis, Birkh\"auser, Boston, 2003.

\bibitem[CoDa93]{CoDa93}
A.~Cohen and I.~Daubechies, {\em Nonseparable bidimensional wavelet bases},
  Rev. Mat. Iberoamericana {\bf 9} (1993), 51--137.

\bibitem[DuSc52]{DuSc52}
R.J. Duffin and A.C. Schaeffer, {\em A class of nonharmonic {F}ourier series},
  Trans. Amer. Math. Soc. {\bf 72} (1952), 341--366.

\bibitem[DuXu01]{DuXu01}
C.F. Dunkl and Y.~Xu, {\em Orthogonal {P}olynomials of {S}everal {V}ariables},
  Encyclopedia of Mathematics and its Applications, vol.~81, Cambridge
  University Press, Cambridge, 2001.

\bibitem[Dut04a]{Dut04b}
D.E. Dutkay, {\em Positive definite maps, representations and frames}, Rev.
  Math. Phys. {\bf 16} (2004), no.~4, 451--477.

\bibitem[Dut04b]{Dut04c}
\bysame
, {\em The local trace function for super-wavelets}, Wavelets,
  Frames, and Operator Theory (Focused Research Group Workshop, College Park,
  Maryland, January 15--21, 2003) (C.~Heil, P.E.T. Jorgensen, and D.~Larson,
  eds.), Contemp. Math., vol. 345, American Mathematical Society, Providence,
  2004, pp.~115--136.

\bibitem[Dut06]{Dut06a}
\bysame
, {\em Low-pass filters and representations of the
  {B}aumslag--{S}olitar group}, Trans. Amer. Math. Soc., to appear,
  http://www.arxiv.org/abs/math.CA/0407344.

\bibitem[Eld02]{Eld02}
Y.C. Eldar, {\em Least-squares inner product shaping}, Linear Algebra Appl.
  {\bf 348} (2002), 153--174.

\bibitem[FJKO05]{FJKO05}
M.~Fickus, B.D. Johnson, K.~Kornelson, and K.A. Okoudjou, {\em Convolutional
  frames and the frame potential}, Appl. Comput. Harmon. Anal. {\bf 19} (2005),
  77--91.

\bibitem[Jor06]{Jor06a}
P.E.T. Jorgensen, {\em Analysis and {P}robability: Wavelets, {S}ignals,
  {F}ractals}, Grad. Texts in Math., vol. 234, Springer-Verlag, New York, to
  appear 2006.

\bibitem[JoKr03]{JoKr03}
P.E.T. Jorgensen and D.W. Kribs, {\em Wavelet representations and {F}ock space
  on positive matrices}, J. Funct. Anal. {\bf 197} (2003), 526--559.

\bibitem[KaRi97]{KaRi97}
R.V. Kadison and J.R. Ringrose, {\em Fundamentals of the {T}heory of {O}perator
  {A}lgebras, {V}ol. {I}: Elementary {T}heory}, Graduate Studies in
  Mathematics, vol.~15, American Mathematical Society, Providence, 1997,
  reprint of the 1983 original Academic Press edition.

\bibitem[KoLa04]{KoLa04}
K.A. Kornelson and D.R. Larson, {\em Rank-one decomposition of operators and
  construction of frames}, Wavelets, Frames and Operator Theory (College Park,
  MD, 2003) (C.~Heil, P.E.T. Jorgensen, and D.R. Larson, eds.), Contemp. Math.,
  vol. 345, American Mathematical Society, Providence, 2004, pp.~203--214.

\bibitem[Nel57]{Nel57}
E.~Nelson, {\em Kernel functions and eigenfunction expansions}, Duke Math. J.
  {\bf 25} (1957), 15--27.

\bibitem[Nel59]{Nel59}
E.~Nelson, {\em Correction to ``{K}ernel functions and eigenfunction
  expansions''}, Duke Math. J. {\bf 26} (1959), 697--698.

\bibitem[PaSc72]{PaSc72}
K.~R. Parthasarathy and K.~Schmidt, {\em Positive {D}efinite {K}ernels,
  {C}ontinuous {T}ensor {P}roducts, and {C}entral {L}imit {T}heorems of
  {P}robability {T}heory}, Lecture Notes in Mathematics, vol. 272,
  Springer-Verlag, Berlin-New York, 1972.

\bibitem[VaWa05]{VaWa05}
R.~Vale and S.~Waldron, {\em Tight frames and their symmetries}, Constr.
  Approx. {\bf 21} (2005), 83--112.

\end{thebibliography}
\ifx\undefined\bysame
\newcommand{\bysame}{\leavevmode\hbox to3em{\hrulefill}\,}
\fi

\end{document}